\documentclass[final,a4paper]{article}
\usepackage{fullpage}
\usepackage[latin1]{inputenc}
\usepackage{amsmath,amsfonts,amsthm}
\usepackage{xspace,tabularx}
\renewcommand{\v}[1]{\boldsymbol{#1}} 
\renewcommand{\i}{\ensuremath{\v{i}}\xspace}
\renewcommand{\j}{\ensuremath{\v{j}}\xspace}
\renewcommand{\k}{\ensuremath{\v{k}}\xspace}
\let\oldmu\mu\renewcommand{\mu}{\v\oldmu}
\let\oldnu\nu\renewcommand{\nu}{\v\oldnu}
  \newcommand{\R}{\ensuremath{\mathbb{R}}\xspace}
\renewcommand{\H}{\ensuremath{\mathbb{H}}\xspace}
  \newcommand{\C}{\ensuremath{\mathbb{C}}\xspace}
\newcommand{\I}{\ensuremath{\boldsymbol{I}}\xspace}
\newtheorem{theorem}{Theorem}
\title{Biquaternion (complexified quaternion) roots of -1}
\author{Stephen J. Sangwine\thanks{The work presented here was carried out at the Laboratoire des
                                   Images et des Signaux, Grenoble, France with financial support
                                   from the Royal Academy of Engineering of the United Kingdom and
                                   the Centre National de la Recherche Scientifique (CNRS).}\\
        Department of Electronic Systems Engineering,\\
        University of Essex, Wivenhoe Park,\\
        Colchester, CO7 9EU, United Kingdom.\\
        Email:~\texttt{S.Sangwine@IEEE.org}}
\date{10 June 2005}
\begin{document}
\maketitle
\begin{abstract}
The roots of -1 in the set of biquaternions (quaternions with complex components, or complex numbers
with quaternion real and imaginary parts) are studied and it is shown that there is an infinite
number of non-trivial complexified quaternion roots (and two degenerate solutions which are the
complex imaginary operator and the set of unit pure real quaternions). The non-trivial roots are
shown to consist of complex numbers with perpendicular pure quaternion real and imaginary parts. The
moduli of the two perpendicular pure quaternions are expressible by a single parameter \textit{via}
a hyperbolic trigonometric identity.
\end{abstract}
\section{Introduction}
In the set of real quaternions \H there are an infinite number of square roots of $-1$, since
$\mu^2=-1$ for any unit pure quaternion $\mu$ \cite[Lemma 1]{arXiv:math.RA/0506034}. This paper
explores the question: what are the roots of $-1$ in the biquaternions, or complexified quaternions,
$\H^\C$. More formally, what is the set of complexified quaternions $\{q\in\H^\C : q^2 = -1 \}$?

Biquaternions, or complexified quaternions, may be represented as quaternions with complex
components, or as complex numbers with quaternion real and imaginary parts, \cite{Ward:1997}:
\begin{align*}
q & = w + \i x + \j y + \k z, && \text{where}\quad w, x, y, z \in \C, \text{and}\quad \i^2=\j^2=\k^2=\i\j\k=-1 \\
  & = q_r + q_i \I,           && \text{where}\quad q_r, q_i   \in \H, \text{and}\quad \I^2=-1\\
\end{align*}
The two representations are equivalent, and conversion between them is straightforward, for example:
$w = S(q_r) + S(q_i)\I$, where $S(\,)$ represents the scalar part of a quaternion and $\I$ is the complex
operator\footnote{We use a capital \I for the complex operator to distinguish it from the first
of the three quaternion operators \i.}; and $q_r = \Re(w) + \Re(x)\i + \Re(y)\j + \Re(z)\k$, where
$\Re(\,)$ is the real part of a complex number. Note that the complex operator \I commutes with the
three quaternion operators $\i,\j,\k$ and therefore with all real quaternions.

We choose to use in this paper the second representation above (a complex number with a real
quaternion real part and a real quaternion imaginary part), and further we represent each of the
real quaternions by a pair of real numbers and a unit pure quaternion:
\[
q = q_r + q_i\I = (a + b \mu) + (c + d \nu)\I, \qquad q_r,q_i\in\H; a,b,c,d\in\R;\mu,\nu\in\H,
S(\mu)=S(\nu)=0,|\mu|=|\nu|=1
\]
This choice makes transparent the solution of the problem at hand, since there are geometric
constraints on the unit pure quaternions, as will be seen, although it hides
the fact that each complexified quaternion has eight real components (each of the two unit pure
quaternions has arbitrary direction in 3-space).

\section{Solution}
\begin{theorem}
Let $q$ be an arbitrary biquaternion (complexified quaternion) represented in the form:
\[
q = (a + b \mu) + (c + d \nu)\I
\]
where $\mu$ and $\nu$ are unit pure real quaternions, \I is the complex imaginary operator 
(therefore $\mu^2=\nu^2=\I^2=-1$), and $a$, $b$, $c$ and $d$ are real. Then the following
solutions, \emph{and no others}, exist for $q^2=-1$:
\begin{align*}
q &= b\mu+d\nu\I,\qquad &&\mu\perp\nu,\quad b^2 - d^2 = 1, \quad a=c=0\\
q &= \mu,               &&\text{that is: }a=c=d=0, b=1\\
q &= \I,                &&\text{that is: }a=b=d=0, c=1\\
\end{align*}
\end{theorem}
The last two solutions are degenerate.
The first solution is the only one in which $q$ is a complexified quaternion.
Note that the constraint on the moduli of $b$ and $d$ in the first solution means
that they are related as hyperbolic sine and cosine by the identity\footnote{Thanks to
Nicolas Le~Bihan for pointing this out.}: $\cosh^2 t - \sinh^2 t = 1$ and therefore the
moduli of the two pure quaternions can be expressed by the single parameter $t$.

The number of possible complexified quaternion solutions to $q^2=-1$ is infinite, since we can
choose $\mu$ anywhere on the unit sphere. There are also an infinite number of choices for $\nu$ in
the plane perpendicular to the chosen $\mu$. We can also choose the moduli of the two unit vectors
in an infinite number of ways to satisfy the constraint $b^2 - d^2 = 1$.

\begin{proof}
The proof first shows that the above solutions are roots of $-1$. We then show that there are
no other solutions by enumeration of the terms obtained when we square $q$. The terms are grouped
into those that are real (which must sum to $-1$), those that are real quaternions (which must sum
to zero) and those that are complex quaternions (which must also sum to zero). Analysis of the
constraints yields the solutions given above and demonstrates that no other solutions exist.

It was shown in \cite[Lemma 1]{arXiv:math.RA/0506034} that any unit pure quaternion is a root of
$-1$, and therefore the unit pure quaternion $\mu$ is a root of $-1$. $\I$ is a root
of $-1$ by definition.

To show that the first solution is a root of $-1$, we square it, remembering that $b$, $d$ and $\I$
commute with $\mu$ and $\nu$, but that $\mu$ and $\nu$ do not commute:
\begin{align*}
(b\mu+d\nu\I)^2 &= b^2\mu^2 + d^2\nu^2\I^2 + bd\mu\nu\I + bd\nu\mu\I\\
                &= -b^2 + d^2 + bd(\mu\nu + \nu\mu)\I
\end{align*}
The term in brackets is zero because the two unit vectors are perpendicular\footnote{The product
of two vectors is the vector product minus the scalar product. The scalar products are zero because
the two vectors are perpendicular, and the dot products cancel.}. The constraint that $b^2 - d^2 = 1$
gives us $q^2=-1$.

We now demonstrate that no other solutions exist. Tabulating the terms obtained when we square $q$
we obtain the following:
\[
\begin{array}{l|rrrr}
q^2      & a          &   b     \mu    &   c \I     &    d     \nu \I \\ \hline
a        & a^2        & a b     \mu    & a c \I     &  a d     \nu \I \\
b \mu    & a b \mu    &  -b^2          & b c \mu \I &  b d \mu \nu \I \\
c     \I & a c     \I & b c     \mu \I &  -c^2      &  c d     \nu \I \\
d \nu \I & a d \nu \I & b d \nu \mu \I & c d \nu \I &    d^2          \\
\end{array}
\]
Note that if $q^2 = -1$:
\begin{enumerate}
\item The four diagonal terms are the only ones that are elements of \R
      and they must sum to $-1$.
\item The two off-diagonal terms in the upper left quadrant ($a b \mu$)
      are the only terms that are elements of \H and since they have the
      same sign, they must be zero. This implies either $a = 0$ or $b\mu = 0$.
      These two cases are considered separately below.
\item All the remaining terms contain \I and must therefore sum to zero.
\end{enumerate}
We now consider the two possibilities $a = 0$ or $b\mu = 0$ which are necessary
but not sufficient conditions to give $q^2 = -1$:
\begin{description}
\item [Case $a=0$:] We retabulate the terms in the table above with the first row and column suppressed:
      \newlength{\tmp}
      \settowidth{\tmp}{$a d \nu \I$}
      \[
      \begin{array}{l|rrrr}
      q^2      &                &   b     \mu    &   c \I     &    d     \nu \I \\ \hline
               &                &                &            &                 \\
      b \mu    &                &  -b^2          & b c \mu \I &  b d \mu \nu \I \\
      c     \I &                & b c     \mu \I &  -c^2      &  c d     \nu \I \\
      d \nu \I & \hspace*{\tmp} & b d \nu \mu \I & c d \nu \I &    d^2          \\
      \end{array}
      \]
      The three diagonal terms must sum to $-1$ and the coefficients of \I must sum to zero.
      This gives us a pair of simultaneous equations, the first in real numbers and the second
      in real quaternions:
      \begin{align*}
      d^2 - b^2 - c^2 &= -1\\
      2 c (b\mu + d\nu) + b d (\mu\nu + \nu\mu) &= 0
      \end{align*}
      Considering the second of these two equations, expressing the products of the two unit
      vectors in terms of dot and cross products we are able to equate scalar and vector parts
      separately to zero:
      \begin{align*}
      - b d (\mu\cdot\nu + \nu\cdot\mu) = -2 b d (\mu\cdot\nu)&= 0\\
      2 c (b\mu + d\nu) + b d (\mu\times\nu + \nu\times\mu) &= 0\\
      \end{align*}
      The scalar part can only be zero if the two unit vectors are perpendicular or if one or both of $b$
      and $d$ is zero. Considering each possibility in turn:
      \begin{description}
      \item [The unit vectors are perpendicular:] In this case, the sum of the cross products will
            cancel in the vector part, leaving us the requirement that $c=0$. This gives us a non-trivial
            root of $-1$, namely:
            \[
            \boxed{q = b\mu+d\nu\I,\quad \mu\perp\nu,\quad b^2 - d^2 = 1}
            \]
      \item [The two unit vectors are not perpendicular:] The scalar part can only be zero if one or
            both of $b$ and $d$ are zero. $b=d=0$ leads only to the degenerate solution $\boxed{q=\I}$.
            If we consider only one of $b$ and $d$ zero, we still require $c=0$ in order for the vector
            part to be zero, otherwise the first term in the second equation does not vanish. Remembering
            that we also have $a=0$, we are left with $q=b\mu$, which if $b=1$, gives the solution
            $\boxed{q=\mu}$; or $q=d\nu\I$, which is not a solution, since it would require $d^2=-1$ which
            is not possible because $d$ is real.
      \end{description}
\item [Case $b\mu = 0$:] We show that this does not lead to any new roots of $-1$ (the only result
      is $q=\I$, which we already have above). We retabulate the terms in the product with
      the second row and column suppressed:
      \settowidth{\tmp}{$b d \nu \mu \I$}
      \[
      \begin{array}{l|rrrr}
      q^2      & a          &                &   c \I     &    d     \nu \I \\ \hline
      a        & a^2        &                & a c \I     &  a d     \nu \I \\
               &            &                &            &                 \\
      c     \I & a c     \I &                &  -c^2      &  c d     \nu \I \\
      d \nu \I & a d \nu \I & \hspace*{\tmp} & c d \nu \I &    d^2          \\
      \end{array}
      \]
      The three diagonal terms must sum to $-1$ and the coefficients of \I must sum to zero.
      Dividing out a factor of two in the second case, and collecting terms we have:
      \begin{align*}
      a^2 - c^2 + d^2 &= -1\\
      a c + (a + c) d \nu &= 0 \implies ac = 0 \quad\text{and}\quad (a + c) d \nu =0
      \end{align*}
      Hence we require: $ac = 0$ and either $a+c=0$ or $d\nu=0$. Considering these possibilities in turn:
      \begin{description}
      \item[Case $a+c=0$:]  This can only be satisfied if $a=0$ and $c=0$, since $a=-c$ would not
                            give $ac=0$. In this case $q$ must be of the form $d\nu\I$, since we
                            already have that $b\mu=0$.
                            The first equation above requires $d^2=-1$ which is impossible,
                            since $d$ is real.
      \item[Case $d\nu=0$:] This can only be satisfied if $a=0$ xor $c=0$, since if both were zero,
                            $q$ would be zero. If $c$ is zero we cannot satisfy $a^2=-1$ since $a$
                            is real. $a=0$ leads to the solution $q = c\I$ which is a degenerate
                            solution if $c=1$, giving us again $\boxed{q=\I}$.
      \end{description}
\end{description}
\end{proof}
\section{Numerical examples}
The non-trivial result given above is not intuitive, so we include three numeric examples.

We can choose $\i$ and $\j$ as examples of unit vectors that are perpendicular, and choose the
modulus of one to be $\sqrt{2}$ and the modulus of the other to be unity so that the difference
of the squared moduli is unity (note that $b$ must be greater than $d$):
$q = \sqrt{2}\i + \j\I$. We then have:
\[
q^2 = (\sqrt{2}\i + \j\I)^2
    = 2\i^2 + \j^2\I^2 + \sqrt{2}(\i\j + \j\i)\I
    = -2 + 1 + \sqrt{2}(\k - \k)\I
    = -1
\]
As a second example, we choose $\mu = \frac{1}{\sqrt{3}}(\i + \j + \k)$ and $\nu =
\frac{1}{\sqrt{2}}(\j-\k)$. To show that these are perpendicular it is only necessary
to compute their dot product $x_{\mu} x_{\nu} + y_{\mu} y_{\nu} + z_{\mu} z_{\nu}$:
\[
\mu\cdot\nu = \frac{1}{\sqrt{3}}\,0 + \frac{1}{\sqrt{3}}\frac{1}{\sqrt{2}} - \frac{1}{\sqrt{3}}\frac{1}{\sqrt{2}} = 0
\]
and it is simple to see that each has unit modulus. We choose\footnote{Other choices are possible
here, but this particular choice makes the algebra simpler.} the moduli to be $\sqrt{3}$ and
$\sqrt{2}$, satisfying $b^2-d^2=1$. Thus $q = (\i + \j + \k) + (\j-\k)\I$. Tabulating the product
terms when $q$ is squared:
\[
\begin{array}{r|rrrrr}
 q^2  &  \i   &   \j &   \k &  \j\I & -\k\I\\ \hline
 \i   &  -1   &   \k &  -\j &  \k\I &  \j\I\\
 \j   & -\k   &  -1  &   \i &   -\I & -\i\I\\
 \k   &  \j   & -\i  &   -1 & -\i\I &    \I\\
 \j\I & -\k\I &  -\I & \i\I &    +1 &  \i  \\
-\k\I & -\j\I & \i\I &   \I & -\i   & +1
\end{array}
\]
Inspection of this table shows that all the terms cancel out, bar a solitary $-1$.

For our third example, we use the same two unit vectors, but in the opposite order, and choose
moduli of $3$ and $2\sqrt{2}=\sqrt{8}$. The difference between these squared moduli is again unity
(note that we could choose two much larger moduli, provided the difference between their squares is
unity). Thus $q = 3\nu + 2\sqrt{2}\mu\I$. Squaring $q$ we obtain:
\[
q^2 = -9 + 8 + 6\sqrt{2}(\nu\mu + \mu\nu)\I
\]
which is $-1$ because the term in the brackets vanishes, $\mu$ and $\nu$ being perpendicular, as is
easily verified by multiplying out $(\j-\k)(\i+\j+\k) + (\i+\j+\k)(\j-\k)$.


\end{document}